\documentclass{jaums}
\usepackage{amscd}
\usepackage{amssymb}
\usepackage{amsmath}

\theoremstyle{thmit} % Numbered and Italic

\theoremstyle{thmrm} % Numbered and Roman

\newtheorem*{oldproof}{Proof}

\newcommand\Calderon{Cal\-der{\'o}n}

\newcommand{\dom}{\operatorname{dom}}
\newcommand{\End}{\operatorname{End}}

\newcommand{\sign}{\operatorname{sign}}
\newcommand{\wt}{\widetilde}

\def\ftoo#1{\stackrel{#1}{\too}}

\def\gS{\Sigma}

\newcommand\gl{\lambda}
\newcommand\gG{\Gamma}
\def\lla{\langle}
\def\noi{\noindent}

\def\rra{\rangle}

\def\too{\longrightarrow}
\def\w{\omega}

\def\wt{\widetilde}

\def\={\cong}
\def\>{\supset}
\def\<{\subset}
\def\ii{^{-1}}
\def\12{\frac{1}{2}}
\def\0{^{\circ}}
\usepackage[english]{babel}
\usepackage[latin1]{inputenc}
\usepackage{times}
\usepackage[T1]{fontenc}
\usepackage{graphicx}

\usepackage{amscd}
\usepackage{amssymb}

\title{Basic Functional Analysis Puzzles of Spectral Flow}
\thanks{{\rm Contribution to the Conference on} {\em Geometry and Quantum Field Theory},
{\rm Max Planck Institute for Mathematics, Bonn, June 20-26, 2010}}

\author{B. Booss-Bavnbek}
\address{IMFUFA, Department of Science, Systems and Models\\Roskilde University\\
Postboks 260, DK-4000 Roskilde, Denmark\\  email: {booss@ruc.dk}
\\ URL: {http://milne.ruc.dk/$\sim$Booss/ }}


\begin{document}

\maketitle

\begin{abstract}
We explain an array of basic functional analysis puzzles on the way
to general spectral flow formulae and indicate a direction of future
topological research for dealing with these puzzles.
\end{abstract}

\noi {\em Dedicated to Alan L. Carey on his 60th birthday}

%% \smallskip

\bigskip

\section{Introduction}\label{s:1}

Over the last decades, substantial progress has been achieved in
analytic approaches to spectral flow in various geometric,
topological and operator algebra settings. For a taste of some
recent results see, {\em e.g.}, M.-T. Benameur {\em et al}.
\cite{BCPRSW}.

Each new approach, each new context displays new and surprising
features, radically new difficulties to be overcome and
astonishing aspects of the new results. How can it be that seemingly
small changes of the setting require different methods and
types of assumptions and yield radically different results?

One explanation can be found in the array of basic functional
analysis puzzles connected with the concept of the spectral flow and
its calculation. In Section \ref{s:2}, we fix the notation and
recall the most elementary spectral flow formula, relating the
symmetric category of curves of self-adjoint Fredholm operators in
separable Hilbert space with the skew-symmetric category of
symplectic functional analysis. Moreover, we shall point to peculiar
functional-analytical properties of geometrically defined operators
like Dirac type operators. We explain which of these properties can
be regained for general elliptic operators and how. In Section
\ref{s:3}, we present our list of basic functional analysis puzzles
on the way to general spectral flow formulae. In Section \ref{s:4} we
indicate a promising further direction to deal with these ``puzzles".

\section{The model case of the functional-analytic approach}
\label{s:2}

To investigate spectral properties of geometrically defined
differential operators like the Laplacian and the Dirac operator on
manifolds with boundary and on partitioned manifolds, one has to
draw on a variety of tools.

\subsection{The von-Neumann approach.} However, some common deep
functional-analytical roots of these formulas have been revealed by
K. Furutani and the author in \cite{BooFur:MIF}, emphasizing the
role of the Cauchy data spaces. More precisely, in the von
Neumann-Kre\u{i}n-Vishik-Birman tradition one is given a complex
separable Hilbert space $\mathcal{H}$ and  a closed symmetric
operator $A$. One defines the symplectic Hilbert space of abstract
boundary values by $\beta(A):=\dom(A^*)/\dom(A)$ with naturally
induced inner product $\lla [x],[y]\rra $ and symplectic form
$\omega([x],[y])$ and the natural {Cauchy data space}
$\operatorname{CD}(A):= \{[x]\mid x\in\ker A^*\}$. One has a
canonical correspondence between all self-adjoint extensions $A_D$
of $A$ with domain $D$ and the Lagrangian subspaces $[D]\<
\beta(A)$. In this framework, {\em e.g.}, if $A_D$ is a
{self-adjoint Fredholm extension} and  $\{C_t\}$ a {continuous
curve} in $\mathcal{B}(\mathcal {H})$ with $\ker(A^*+C_t+s)\cap
\dom(A)=\{0\}$ for small $|s|$ ({weak inner UCP}), one obtains that
$\{\operatorname{CD}(A+C_t),[D]\}$ is a continuous curve of Fredholm
pairs of Lagrangians and $\operatorname{SF}\{(A+C_t)_D\} =
\operatorname{MAS}\{\operatorname{CD}(A+C_t),[D]\}$, relating the
spectral flow of a self-adjoint Fredholm operator under bounded
variation with the Maslov index of the corresponding curve of
Lagrangians in the abstract boundary space.

The strength of this functional-analytical approach turns up when
dealing with systems of ordinary differential equations on the
interval, generalizing the classical Morse index theorem for
geodesics on Riemannian manifolds to Subriemannian manifolds. It
recovers  the Floer-Yoshida-Nicolaescu splitting results for the
spectral flow of curves of Dirac operators on partitioned manifolds
(i.e., the family version of the Bojarski Conjecture), and it
provides a basic functional-analytical model for {\em quantization}
and {\em tunneling}, relating {\em spectral} and {\em symplectic}
invariants.

\subsection{The challenge of varying domain.}
Then, how can we transgress the limitations of the general
functional-analytical approach? What if we don't keep
the domain fixed under variation; nor restrict to bounded (i.e., 0
order) perturbations; nor confine the applicability to ordinary
differential
 equations or Dirac   type operators with constant coefficients in normal direction
  (product case) close to boundary? A series of recent papers took up the challenge of
perturbed Dirac operators and general linear elliptic differential
operators; investigated weak inner UCP; established the existence of
self-adjoint Fredholm extensions; admitted variation of domain and
skew boundaries; and investigated uniform structures and continuous
deterministic and random perturbations, see, {\em e.g.}, A. Axelsson
{\em et al}. \cite{AKM:QEF}; joint work of the author with G. Chen,
M. Lesch and C.~Zhu in \cite{CheZhu:PPS}, \cite{BoLe09},
\cite{BoLeZh09}, \cite{BooZhu:MI}; J. Eichhorn \cite{Eic06}; F.
Gesztesy {\em et al}. \cite{GLMST}; and J. Sj{\"o}strand
\cite{Sjo:EDN}. We conclude that the ``natural" (von Neumann)
approach is insufficient, and more {analysis} ({\em e.g.}, splitting
the coefficients near the boundary and {pseudodifferential}
calculus) is needed.

\subsection{Seeley's \Calderon\ projection and Dirac operator folklore.}
Let $M$ be a smooth compact Riemannian manifold with boundary $\gS$,
$E,F$ Hermitian vector bundles over $M$, and $A:C^\infty(M,E)\to
C^\infty(M,F)$ an elliptic differential operator (of first order).
Recall that $\rho:L^2_s(M,E)\to L^2_{s-1/2}(\gS,E|_{\gS})$ for $s>
1/2$ is extendable to ${\mathcal D}_{\operatorname{max}}(A)$. Then
the classical definition of the {Cauchy data space} $N_+^0(A)$ of
$A$ is the closure of$\{\rho u\mid Au=0 \text{ in } M\setminus\gS,\,
u\in C^\infty(M,E)\}$ in ${L^2(\gS,E|_{\gS})}$. R.T. Seeley
\cite{See:SIB,See:TPO} proved that this Cauchy data space can be
obtained as the range of a pseudodifferential projection. The basic
ingredients for Seeley's result have been the construction of an
invertible extension $\wt A$ of $A$ over a closed manifold $\wt M$
by extending $A$ to a collar, then doubling and applying symbolic
calculus and UCP management. As a result, he received a Poisson
operator $K_\pm:=\pm r^\pm \wt A\ii\rho^* J(0)$ where
$J(0)=\sigma(A)(\cdot,\nu)\in \End(E|_{\gS})$ denotes the principal
symbol of $A$ in normal direction at the boundary. He showed that
the operator $C_\pm := \rho K_\pm$ is a pseudodifferential
projection onto $N^0_+(A)$ and called it the {Calder{\'o}n
projection}.

It was shown by K.P. Wojciechowski and the author \cite[Chapters 9
and 12]{BooWoj:EBP} that Seeley's construction is canonical (i.e., {\em natural},
{\em explicit}, {\em transparent}, and {\em free of choices}) for
Dirac type operators when the metric structures are product close to
the boundary. As a consequence, we obtained the Lagrangian property
of the Cauchy data space. The reason is that for such operators the
invertible extension $\wt A$ can be explicitly defined on the very
closed double $\wt M$ of $M$ - without inserting additional collar
near the boundary and not involving any other choices. As a
consequence, the Cauchy data spaces, respectively, the \Calderon\
projection varies continuously under smooth deformation of the data
defining the Dirac operator, proved by M. Lesch, J. Phillips, and
the author in \cite{BooLesPhi:UFOSF}.

These results for Dirac type operators can be traced back to the
``{Dirac Operator Folklore}": (i) weak inner UCP, i.e., $\ker A\cap
{\dom}(A_{\operatorname{min}})=\{0\}$ with
${\dom}(A_{\operatorname{min}})=L^2_{1,\operatorname{comp}}(M,E)$;
(ii) symmetric principal symbol of the tangential operator $B$ in
the decomposition $A=J_0(\partial_x+B)$ where $x$ denotes the inner
normal variable; and (iii) a precise invertible double. From that
alone, one can derive the transparent definition of the
{Calder{\'o}n projection}, the Lagrangian property of the Cauchy
data space, the existence of a self-adjoint Fredholm extension given
by a regular pseudodifferential boundary condition, the Cobordism
Theorem, and the continuous dependence of input data.

One may wonder, how special are operators of Dirac type compared to
arbitrary linear first order elliptic differential operators? The
short answer is that property (i) may be lost but is indispensable,
hence must be assumed. Property (ii) implies property (i) (if it is
valid for the tangential operators on arbitrary hypersurfaces), but
else it is dispensable (for details see below Section \ref{s:3}).
Property (iii) can be maintained by replacing Seeley's classical
construction by a new construction, inspired by B. Himpel {\em et
al}. \cite{HimKirLes:CPH} and worked out in \cite{BoLeZh09}.

\subsection{The invertible double, revisited}
We summarize the new construction. First, we bring a given general
elliptic differential operator of first order in product form $A =
{J\bigl(\partial_{x}+B\bigr)}$ close to the boundary by suitable
choice of the metric. Here, $J$ and $B$ vary with the normal
variable $x$. Note that dropping the geometric Dirac operator
context, the metric structures need no longer to be fixed.

We obtain a canonical new invertible double $\tilde A_{T}$ with
\[
\dom(\tilde A_{T}):=\{\binom{e}{f}\in L^2_1(M,E\oplus F)\mid \varrho
f=T\varrho e\}, \] where $\tilde{A} : C^\infty(M,E\oplus
F)\ftoo{A\oplus(-A^t)} C^\infty(M,F\oplus E)$ and $T
\in\operatorname{Hom}(\gS,E|_\gS,F|_\gS)$ invertible bundle
homomorphism with  $J_0^*T$ positive definite. Then $\tilde{A}_{T}$
is a {Fredholm operator with compact resolvent} with
$\ker\tilde{A}_{T} = Z_{+,0} \oplus Z_{-,0}$ and
$\operatorname{coker} \tilde{A}_{T} \simeq Z_{-,0} \oplus Z_{+,0}$
where $Z_{+,0} := \{f\in L^2_1(M,E)\mid Af=0, \varrho f=0\}$ and
$Z_{-,0}$ denotes the corresponding kernel of $A^t$\,. For the most
part of our work we pick $T := (J_0^t)^{-1}$\,. Denoting the
pseudo--inverse of $\tilde{A}_{T}$ by  $\tilde{G}$ , we define
Poisson operators $K_{\pm}  :=  \pm\, r^{\pm}\tilde{G}\varrho^*J_0:
L^2_s(\Sigma,E)\to L^2_{s+\12}(M,E)$ $(L^2_{s+\12}(M,F))$ and
\Calderon\ operators $C_{+} := \varrho_+K_+, \quad C_{-}  :=
T^{-1}\varrho_{-}K_{-}$. We obtain that $C_{\pm}$ are projections
with $C_+ + C_{-} = I$ and $C_+(L^2)  =  N_+^0,\  C_-(L^2)   =
T^{-1} N_-^0$\,.

The most delicate part of the new construction is the investigation
of the mapping properties of the pseudo-inverse $\wt G$, the Poisson
operators $K_\pm$ and the \Calderon\ projection $C_\pm$\,.

Our model operator is $A=J\bigl(\frac{\partial }{\partial x}
+B(x))+0.\ \text{order}$. From the ellipticity of $A$ we have that
$i\xi+B(x)$ is invertible for real $\xi$ of sufficiently large
numerical value (ray of minimal growth). We put
$Q_+(x):=\frac{1}{2\pi i} \int_{\Gamma_+} e^{-x\lambda} (\lambda-
B(0))^{-1}d\lambda$ a family of sectorial projections where
$\Gamma_+$ is an infinite contour which encircles the eigenvalues of
$B(0)$ in the right half plane. We notice that $Q_+(x)$ corresponds
to $e^{-x B(0)}1_{[0,\infty)}(B(0))$ if $B(0)=B(0)^*$. We had to
display a delicate balance on a knife edge between general operator
theory and pseudodifferential calculus when we realized that {\em a
priori} $Q_+(x)=O(\log x), x\to 0+$, hence $P_+:=Q_+(0)$ is possibly
unbounded. Within the pseudodifferential calculus, it follows,
however, from T. Burak \cite{Bur:SPE}, K.P. Wojciechowski
\cite{Woj85}, V. Naza{\u\i}kinski{\u\i\ {\em et al.}
\cite{NSSS:SPV}, and R. Ponge \cite{Pon:SAZ} (with minor, but
necessary additions and corrections in \cite{CheZhu:PPS}) that
\[
P_+\, :=\, Q_+(0) \, =\, \frac{-1}{2\pi i}B(0)\int_{\gG_+}\,\gl\ii\,
(B(0)-\gl)\ii\, d\gl
\]
is a bounded pseudodifferential projection. {\em A posteriori}, we
obtain $Q_+(x)\to P_+$ strong, $x\to 0+$.

Another hopefully useful concept introduced in \cite{BoLeZh09} is
{the approximative Poisson operator}
$R:C^\infty(\gS,E|_{\gS})\longrightarrow C^\infty(\mathbb{R}_+\times
\gS,E|_{\gS}\oplus F|_{\gS})$ with $R\xi
(x):=\varphi(x){Q_+(x)\xi\choose T Q_-(x)\xi}$, where $\varphi$ is a
suitable { cut-off function at } 0. One finds $R=\tilde
A_T^{-1}\varrho^*$ +{ regularising remainder}. That permits to
analyze the mapping property of $R:L_s^2(\gS,E|_{\gS})\to
L^2_{s'}(\mathbb{R}_+\times \gS,E|_{\gS}\oplus F|_{\gS})$ in
dependence of $A$.

\subsection{A recent result on sectorial projections.}
Regarding uniform structures, it turns out that $C_+(A)- P_+(B(0))$
is a pseudo\-differential operator of order $-1$ and that $A\mapsto
C_+(A)$ is as regular as $A\mapsto P_+(B(0))$ under the condition
$\dim Z_0(A), \dim Z_0(A^t)=\text{const}$. Now, \cite[Theorem
1.1]{CheZhu:PPS} proves that the sectorial projections of elliptic
semi-classical pseudodifferential operators on closed manifolds
depend continuously on the initial operator in a natural Fr{\'e}chet
topology, if there exist suitable spectral cuttings for the
principal symbol (like no purely imaginary eigenvalues of the
principal symbol, which is exactly satisfied for the tangential
operator $B(0)$ of any elliptic operator $A$ over a compact smooth
manifold with boundary). Consequently, the \Calderon\ projection
$C_+(A)$ varies continuously in the operator norm of bounded
operators on $L^2(\gS,E|_{\gS})$, if the coefficients of $A$ and all
its derivatives vary continuously. Moreover one obtains that
$(A,P)\mapsto A_P$ is continuous in graph topology, if $P$ runs in
the space of ``regular" boundary conditions.

Further applications for $A=A^t$ are that the Cauchy data space is
Lagrangian in the Hermitian symplectic Hilbert space
$L^2(\gS,E|_{\gS}), \langle\cdot,J(0)\cdot\rangle)$; the existence
of a self-adjoint Fredholm extension $A_{C_\pm}$ (for a suitable
choice of the auxiliary bundle homomorphism $T$ over $\gS$); and the
cobordism invariance of the index for arbitrary symmetric elliptic
differential operators on closed manifolds: $\sign iJ(0)$ vanishes
on $\bigoplus\limits_{\lambda \text{ imaginary}}\ker
(B(0)-\lambda)^N$, $N\gg 0$.

\section{A personal choice of functional analysis puzzles}
\label{s:3}

From the preceding summary we can extract an array of
functional-analytical puzzles on the way to general spectral flow
formulae.

\subsection{Geometrically defined vs. general coefficients.} In
applications, there is a decisive difference between {\em ad-hoc}
models and models based on first principles, as pointed out, {\em
e.g.}, by Y. Manin \cite{Man:MKI}. {\em Ad-hoc} models are based on
fancied hypotheses about the interrelation between different
features and on estimates of the rates and other coefficients.
Mathematically speaking, they require general coefficients. On the
contrary, equations and coefficients in theoretically based models
have a direct meaning, {\em e.g.}, when derived from minimal
principles. Often, to exploit this meaning one better restricts the
consideration to geometrically defined operators, instead of
striving for the goal of ``highest generality". Clearly, for gaining
mathematical insight both approaches have their merits and yield
their own way of {\it transparency}. In the example presented above
in Section \ref{s:2}, the Dirac case yields a simple construction of
the invertible double while the general approach yields a list of
universal essentials for getting through.

\subsection{Fixed operator vs. deformation curve.} Addressing curves
instead of single points is as old as celestial mechanics and
variational calculus. To embed such questions in a systematic way
into a family setting of deformations is not a new idea; it goes
back to J.L. Lagrange's second letter (in Latin) to Euler regarding
the derivation of what is called now the {\em Euler-Lagrange
Equation} \cite{Lag55}. Following Lagrange, it seems a tenet of the
mathematics of our time to address deformation questions at the
first place. As a rule it turned out, {\em e.g.}, in Index Theory
that family versions are more demanding than single operator
formulae. In contractible spaces the situation is different when,
{\em e.g.}, the spectral flow of a curve solely depends on the
endpoints. Then, like in Lagrange's idea, the embedding of a problem
into a deformation curve may facilitate the treatment and not
complicate.

\subsection{Bounded vs. unbounded operators.} With some right, we may
forget about that distinction when working with an  elliptic
operator $A$ (say symmetric and of order 1) on a closed manifold
$M$. Then there is no difference between minimal and maximal domain.
It is always equal to the Sobolev space $L_1^2(M)$. Moreover, in
that case the Riesz transform $A\mapsto A(I+A^2)^{-1/2}$ yields a
bounded operator in $L^2(M)$ and is continuous in suitable operator
norms, see \cite[Chapter 17]{BooWoj:EBP}. The situation is much more
blurred for elliptic operators on manifolds with boundary. There,
the general functional analysis picture has strongly
counter-intuitive traits.

Let $\mathcal{CF}(H)$ denote the space of closed (not necessarily
bounded) Fredholm operators in a fixed complex separable Hilbert
space $H$ and let $\mathcal{CF}^{\operatorname{sa}}(H)$ denote the
subspace of self-adjoint elements. For index theory, H.O. Cordes and
J.P. Labrousse \cite{CorLab} have shown that the index is constant
on the connected components of $\mathcal{CF}(H)$ and yields a
bijection between the integers and the connected components. For the
spectral flow, quite a different result was proved in
\cite{BooLesPhi:UFOSF}: While the space of bounded self-adjoint
Fredholm operators decomposes in three connected components (the
contractible spaces of essentially positive, respectively
essentially negative operators and the non-trivial component with
homotopy type of Bott periodicity), the space
$\mathcal{CF}^{\operatorname{sa}}(H)$ is connected and its homotopy
type is not fully revealed. Moreover, equipping the space
$\mathcal{CF}^{\operatorname{sa}}(H)$ with the graph (gap) topology
and the space of bounded operators with the operator norm, the Riesz
transform is not continuous as shown by a counterexample provided by
B. Fuglede (for details see {\em l.c.}).

\subsection{Self-adjoint vs. general.} Motivated by the method of
replacing a differential equation by difference equations, D.
Hilbert and R. Courant \cite{CoHi} expected ``linear problems of
mathematical physics which are correctly posed to behave like a
system of $N$ linear algebraic equations in $N$ unknowns... If for a
correctly posed problem in linear differential equations the
corresponding homogeneous problem possesses only the trivial
solution zero, then a uniquely determined solution of the general
inhomogeneous system exists. However, if the homogeneous problem has
a nontrivial solution, the solvability of the inhomogeneous system
requires the fulfillment of certain additional conditions." This is
the {\em heuristic principle} which Hilbert and Courant saw in the
{\em Fredholm Alternative}. G. Hellwig in the real setting and I. N.
Vekua in complex setting (both nicely explained in the recent H.
Kalf \cite{Ka}) disproved it in 1952. Independently of each other they
discovered symmetric differential operators on the disc with
non-self-adjoint boundary condition where the Fredholm Alternative
fails.

From the chiral splitting of Dirac type operators we have learnt
that self-adjoint and non-self-adjoint problems can be
related to each other. One instant is the Cobordism Theorem for two
linear elliptic, not necessarily symmetric operators on closed
manifolds which appear as components of the tangential operator for
a self-adjoint boundary problem, \cite[Corollary 21.6]{BooWoj:EBP}.

It is remarkable how easy it is to apply the Spectral Theorem to
prove the continuous dependence of spectral projections outside a
spectral cut for {\em symmetric} elliptic differential operators on
closed manifolds (see \cite[Proposition 7.15]{BoLeZh09}) and how
elaborate the arguments become for proving a similar result {\em
without} symmetry assumptions (see \cite{CheZhu:PPS}).

\subsection{Functional analysis vs. pseudodifferential analysis.} The
investigation of the mapping properties for constructing sectorial
and \Calderon\ projections from elliptic operators yields a treasure
of situations where claims can be formulated in general
functional-analytical terms but be proved only by advanced
pseudodifferential analysis. As examples, see the preceding
discussion of the boundedness of the sectorial projection
$P_+(B(0)$); the coincidence of the mapping property of $B(0)\mapsto
P_+(B(0))$ and $A\mapsto C_+(A)$; and the mentioned recent delicate
proof of the continuous dependence of $P_+(B(0))$ on $B(0)$.

\subsection{Strong symplectic vs. weak symplectic.} From classical
mechanics and the usual treatment of Dirac operators, we are
accustomed to strong symplectic structures, i.e., we assume that the
symplectic form $\omega$ can be written as a scalar product
$\w(x,y)=\lla Jx,y\rra$ with bounded invertible (i.e., also the
inverse is bounded) generator operator $J$. On a smooth compact
manifold $M$ with boundary $\gS$, any elliptic operator $A$ (say of
order 1 and symmetric) induces strong symplectic  structures on the
von-Neumann boundary value space $\beta$ defined above in Section
\ref{s:2} and on $L^2(\gS)$ with $J$ defined by the principal symbol
of $A$ over $\gS$ in inner normal direction. Formally in the same
way, we obtain a symplectic structure for the Sobolev space
$L^2_{1/2}(\gS)$ where all the boundary values of the domain of the
extensions of $A$ are placed by Sobolev restriction. However, for
$\dim\gS\ge 1$, that structure is no longer strong but becomes weak,
see \cite[Section 2, Remark]{BooZhu:MI}. In weak symplectic
analysis, we don't know whether the space of Lagrangian subspaces is
contractible; whether the homotopy of the space of Fredholm pairs of
Lagrangian subspaces is of Bott periodicity; nor whether there exist
Fredholm pairs of Lagrangian subspaces with negative index, see
\cite[Section 2.3]{BooZhu:MI}.

\subsection{Weak inner UCP?} For operators of Dirac type, the weak
Unique Continuation Property can be obtained in two different ways,
either by exploiting that the principal symbol of the Dirac
Laplacian is in diagonal form and real or by exploiting that the
principal symbol of the tangential operators are symmetric for all
hypersurfaces, see \cite[Chapter 8]{BooWoj:EBP} for details or
\cite[Theorem 1.3]{BoLe09} for outlines and references. In
difference to the usual Unique Continuation Property for elements
belonging to the kernel of an elliptic operator, the property {\em
weak inner UCP}, discussed above in Section \ref{s:2} is purely
functional-analytical. As an immediate consequence, C. Zhu and the
author obtained the local stability of weak inner UCP, for the
references and wider ramifications see \cite[Section 4]{BoLe09}. The
stability of weak (global) UCP was obtained by the author and M.
Marcolli and B. Wang \cite{BooMarWan:WUCP} for {\em mild} non-linear
perturbations of the Dirac operator, motivated by Seiberg--Witten
Theory.

\bigskip

We shall not elaborate on the many other puzzles. For instance, one
may wonder about the functional-analytical roots of the noted
differences between {\em homotopy invariance}, valid for index and
spectral flow in suitable setting, and solely {\em spectral
invariance} of $\eta$-invariant and $\zeta$-function regularized
determinants. Another puzzle, not addressed here, are the
differences and relations between the desuspension character of
spectral flow formulae going a dimension down (mostly rather
delicate from an analysis point of view) and the suspension
character of rather different spectral flow formulae, going a
dimension up (and often more easily accessible). Since the first
tries by K.P. Wojciechowski and the author in the early 1980's
(quoted in \cite[Theorem 17.13 vs. Theorem 17.17]{BooWoj:EBP}) these
questions have been studied extensively for suspended actions. In
particular, I refer to the programmatic V. Mathai \cite{Mat:SFE} and
the follow-up papers, e.g., by N. Keswani \cite{Kes:NEI} and the
recent M.-T. Benameur and P. Piazza \cite{BenPia:IER}.

\section{How to deal with these puzzles?}\label{s:4}

History of mathematics (and of sciences, as well) provides ample
evidence of changes between periods of expansion (diversification) and
periods of consolidation (establishing deep, principal interrelations).
A famous case is, how the ideas of R. Bott, F. Hirzebruch, I.M. Singer,
and M.F. Atiyah (and followers) lead to the identification of
{\em Fredholm operators} and {\em index problems} in wide fields of geometry
and a corresponding unprecedented interconnection between topology, geometry,
functional analysis, PDEe, dynamical systems, number theory, and
mathematical physics. Similarly, one may expect that the avalanche of new results
on spectral invariants of operator curves, though pointing in many seemingly
unrelated directions, will help to single out one or two key concepts for
dealing with the listed (and supplementary) ``puzzles" around spectral flow.

To overcome - or better to make maximal use of - the vaste amount of
inspiring, but spread  calculations, it will not be easy to single
out a specific direction of dealing with all the puzzles in one
round. One candidate for such a unifying approach is the
concentration on the homotopy type of the operator spaces involved.
For K.P. Wojciechowski and me, that was the starting point of our
joint work, see {\em e.g.} \cite[Chapters 15-17]{BooWoj:EBP}. The
task is easy to formulate: look for the involved subspaces of
unitary operators and check whether Bott periodicity is maintained,
respectively determine deviations in homotopy type, and do it both
in general functional analysis terms and in pseudodifferential
operator terms. To me, the work, {\em e.g.}, by P. Kirk and M. Lesch
\cite[Sections 2 and 6]{KiLe00} indicates that this program
continues to be promising.

% The Acknowledgements are an un-numbered section
\section*{Acknowledgements}
% Acknowledgements text here
For the preceding considerations I must take the responsibility
alone, assuming that most of the views were shared or will be shared
by Alan Carey. The findings are also based on continuing discussions
with students and colleagues. In particular, I am indebted to Kenro
Furutani (Tokyo), Matthias Lesch (Bonn), Ryszard Nest (Copenhagen),
John Phillips (Victoria), Chaofeng Zhu (Tianjin) and my late
collaborator for decades Krzysztof Wojciechowski. I thank also the
referees for helpful suggestions.

%% \end{document}

\end{document}